\documentclass[oneside,11pt]{amsart}
\usepackage{amsmath, amsfonts,amsthm,times,graphics}
\usepackage{mathrsfs}
\usepackage[usenames]{color}
\usepackage[active]{srcltx}
 \makeatletter
\renewcommand*\subjclass[2][2000]{%
  \def\@subjclass{#2}%
  \@ifundefined{subjclassname@#1}{%
    \ClassWarning{\@classname}{Unknown edition (#1) of Mathematics
      Subject Classification; using '1991'.}%
  }{%
    \@xp\let\@xp\subjclassname\csname subjclassname@#1\endcsname
  }%
}
 \makeatother

\newtheorem*{ThmA}{Theorem A}

\newtheorem{Thm}{Theorem}[section]
\newtheorem{Cor}[Thm]{Corollary}

\theoremstyle{definition}

\theoremstyle{remark}

\numberwithin{equation}{section}

\newcommand{\ee}{\mathrm{e}}

\theoremstyle{definition}

\def\be{\begin{equation}}
\def\ee{\end{equation}}

\newcommand{\ben}{\begin{enumerate}}
\newcommand{\een}{\end{enumerate}}

\newcommand{\br}{\begin{rem}}
\newcommand{\er}{\end{rem}}
\newcommand{\brs}{\begin{rems}}
\newcommand{\ers}{\end{rems}}
\newcommand{\bo}{\begin{obser}}
\newcommand{\eo}{\end{obser}}
\newcommand{\bos}{\begin{obsers}}
\newcommand{\eos}{\end{obsers}}
\newcommand{\bpf}{\begin{pf}}
\newcommand{\epf}{\end{pf}}
\newcommand{\ba}{\begin{array}}
\newcommand{\ea}{\end{array}}
\newcommand{\beq}{\begin{eqnarray}}
\newcommand{\beqq}{\begin{eqnarray*}}
\newcommand{\eeq}{\end{eqnarray}}
\newcommand{\eeqq}{\end{eqnarray*}}

\numberwithin{equation}{section}

\newcounter{minutes}
\divide\time by 60
\newcounter{hours}
\multiply\time by 60 \addtocounter{minutes}{-\time}
\begin{document}

\bibliographystyle{amsplain}
\title [] {A Hardy-Littlewood type Theorem and a Heinz type inequality}

\author{Shaolin Chen}
\address{S. L. Chen,    College of Mathematics and
Statistics, Hengyang Normal University, Hengyang, Hunan 421008,
People's Republic of China} \email{mathechen@126.com}

\author{Hidetaka Hamada$^{\ast}$}
\thanks{$^{\ast}$ Corresponding author}
\address{H. Hamada, Faculty of Science and Engineering, Kyushu Sangyo University,
3-1 Matsukadai 2-Chome, Higashi-ku, Fukuoka 813-8503, Japan.}
\email{ h.hamada@ip.kyusan-u.ac.jp}

\author{Dou Xie}
\address{D. Xie,    College of Mathematics and
Statistics, Hengyang Normal University, Hengyang, Hunan 421008,
People's Republic of China} \email{2393556505@qq.com}


\date{\today  
File: Lips. spaces. tex}

\begin{abstract}
The main aim of this paper is to investigate the Hardy-Littlewood type Theorem and
the Heinz type inequality on functions induced by a differential operator.
We first  prove a more general Hardy-Littlewood type theorem for the Dirichlet solution of a differential operator
which depends on $\alpha >0$
over the unit ball $\mathbb{B}^n$ of $\mathbb{R}^n$ with $n\geq 2$,
related to the Lipschitz type space defined by a fast majorant.
We find that the case $\alpha>0$ is completely different from the case $\alpha=0$.
Then a more general Heinz type inequality for the Dirichlet solution of a differential operator will also be established
in the case $\alpha>n-2$.
\end{abstract}






\maketitle \pagestyle{myheadings} \markboth{S. L. Chen, H.
Hamada and D. Xie}{A Hardy-Littlewood type Theorem and a Heinz type inequality}

\section{Introduction}\label{sec1}
For an  integer $n\geq2$, let $\mathbb{R}$ and $\mathbb{R}^{n}$ be
the set of real numbers
 and the usual real vector space of dimension $n$, respectively. 
For $y=(y_{1},\ldots,y_{n})$, $x=(x_{1},\ldots,x_{n})\in\mathbb{R}^{n}$, we define
the Euclidean inner product $\langle \cdot ,\cdot \rangle$ by
$\langle x,y\rangle=\sum_{k=1}^{n}x_{k}y_{k}$ so that the Euclidean
length of $x$ is defined by
$$|x|=\langle x,x\rangle^{1/2}=\left(\sum_{k=1}^{n}|x_{k}|^{2}\right)^{1/2}.
$$
Denote a ball in $\mathbb{R}^{n}$ with center
$x_{0}\in\mathbb{R}^{n}$ and radius $r$ by
$\mathbb{B}^{n}(x_{0},r)$. In particular, let
$\mathbb{B}^{n}:=\mathbb{B}^{n}(0,1)$,
$\mathbb{S}^{n-1}:=\partial\mathbb{B}^{n}$ and $\mathbb{D}:=\mathbb{B}^{2}$.

For $x\in\mathbb{B}^{n}$, let $$\delta_{\alpha}(x)=(1-|x|^{2})^{\alpha}$$
be the standard weight in $\mathbb{B}^{n}$, where $\alpha\in\mathbb{R}$. Let us recall the
differential operator $\Delta_{\alpha}$ related to the  standard weight $\delta_{\alpha}$ as follows:
\be\label{C-y}\Delta_{\alpha}:={\rm div}(\delta_{\alpha}^{-1}\nabla)+\alpha(n-2-\alpha)\delta_{\alpha+1}^{-1},\ee
where the symbols $\nabla$ and ${\rm div}$ denote the gradient and divergence, respectively (cf. \cite{A-P,B-H,LP,LP-2009,LPS,O-1}). For $n=2$,
the formula (\ref{C-y}) points at a relation to the so-called conductivity equations considered by
Astala and P\"aiv\"arinta \cite{A-P} (see also \cite{O-1}).
For a function $\varphi$ on $\mathbb{S}^{n-1}$, we denote by $P_{\alpha}[\varphi]$ the Dirichlet solution, for the  differential operator $\Delta_{\alpha}$,
 of $\varphi$ over $\mathbb{B}^{n}$, that is,   $\Delta_{\alpha}P_{\alpha}[\varphi]=0$  in $\mathbb{B}^{n}$ and $P_{\alpha}[\varphi]=\varphi$ on $\mathbb{S}^{n-1}$, where $\alpha>-1$.
For $\varphi\in C(\mathbb{S}^{n-1})$, it follows from  \cite[Theorem 2.4]{LP} (or \cite[Lemma 2.9]{LPS}) that

$$P_{\alpha}[\varphi](x)=\int_{\mathbb{S}^{n-1}}\mathbf{P}_{\alpha}(x,\zeta)\varphi(\zeta)d\sigma(\zeta),~x\in\mathbb{B}^{n},$$
where
$\sigma$ is the rotation-invariant positive Borel measure on $\mathbb{S}^{n-1}$ for which $\sigma(\mathbb{S}^{n-1})=1$,
$$\mathbf{P}_{\alpha}(x,\zeta)=C_{\alpha}\frac{(1-|x|^{2})^{1+\alpha}}{|x-\zeta|^{n+\alpha}}$$ is the Poisson type kernel, $\alpha>-1$ is a  real parameter
and
$$C_{\alpha}=\frac{\Gamma\left(\frac{n+\alpha}{2}\right)\Gamma\left(1+\frac{\alpha}{2}\right)}{\Gamma\left(\frac{n}{2}\right)\Gamma(1+\alpha)}.$$
Here
$\Gamma(s)=\int_{0}^{\infty}t^{s-1}e^{-t}dt$  for $\mbox{Re}(s)>0$ is the standard Gamma function.
By the analytic continuation, the Gamma function $\Gamma$ has no zeros
and is a meromorphic function
which is holomorphic in the whole complex plane except zero and the negative integers,
where this function has simple poles.

Note that $P_{\alpha}[\varphi]$ is continuous on $\overline{\mathbb{B}^n}$, where  $\varphi\in C(\mathbb{S}^{n-1})$ (see \cite{LP,LP-2009,LPS}).
In particular, if $\alpha=0$, then $P_{0}[\varphi](x)$ is
 harmonic for $x\in\mathbb{B}^{n}$.
It  has also attracted much attention of many authors
\cite{B-H,O-1,O-2} on the
differential operator (\ref{C-y}) for $n=2$.

Throughout of this paper,
 we use the symbol $M$ to denote the various positive
constants, whose value may change from one occurrence to another.

\section{Preliminaries and main results}\label{sec2}

A continuous increasing function $\omega:[0,\infty)\rightarrow[0,\infty)$ with $\omega(0)=0$ is called a majorant if
$\omega(t)/t$ is non-increasing for $t\in(0,\infty)$ (see \cite{Dy1,Dy2}).
Then for any $c\geq 1$ and $t>0$, we have
$\omega(ct)\leq c\omega(t)$.
For $\lambda_{0}>0$ and $0<\lambda<\lambda_{0}$, we consider the
following conditions on a majorant $\omega$:
\be\label{eq2x}
\int_{0}^{\lambda}\frac{\omega(t)}{t}\,dt\leq\, M\omega(\lambda)
\ee
and
\be\label{eq3x}
\lambda\int_{\lambda}^{\infty}\frac{\omega(t)}{t^{2}}\,dt\leq\, M
\omega(\lambda),
\ee
where $M$ denotes a positive constant. A majorant $\omega$ is henceforth  called fast (resp. slow) if condition (\ref{eq2x}) (resp. (\ref{eq3x}) ) is fulfilled.
In particular, a majorant $\omega$ is said to be  regular if it satisfies the
conditions (\ref{eq2x}) and (\ref{eq3x}) (see \cite{Dy1,Dy2}).

Given a majorant $\omega$ and a subset  $\Omega$ of
$\mathbb{R}^{n}$, a function $f$ of  $\Omega$ into $\mathbb{R}$ is
said to belong to the  Lipschitz  space
$\mathscr{L}_{\omega}(\Omega)$
 if there is a positive constant $M$ such that

 \be\label{rrt-1}|f(x)-f(y)| \leq\,M\omega\left(|x-y|\right), \quad
x,y \in \Omega.\ee
Furthermore, let
$$\|f\|_{\mathscr{L}_{\omega}(\Omega),s}:=\sup_{x,y\in\Omega,x\neq\,y}\frac{|f(x)-f(y)|}{\omega(|x-y|)}<\infty.$$
Note that if $\Omega$ is a proper subdomain of $\mathbb{R}^{n}$ and
$f\in \mathscr{L}_{\omega}(\Omega)$, then $f$ is continuous on
$\overline{\Omega}$ and (\ref{rrt-1}) holds for $x, y \in
\overline{\Omega}$ (see \cite{Dy2}).  In particular, we say that
a function $f$ belongs to the  local Lipschitz space
$\mbox{loc}\mathscr{L}_{\omega}(\Omega)$ if (\ref{rrt-1})
holds, with a fixed
positive constant $M$, whenever $x\in \Omega$ and
$|x-y|<\frac{1}{2}d_{\Omega}(x)$ (cf. \cite{Dy2,GM,L}). Moreover, $\Omega$ is called an  $\mathscr{L}_{\omega}$-extension domain if
$\mathscr{L}_{\omega}(\Omega)=\mbox{loc}\mathscr{L}_{\omega}(\Omega).$ On the geometric
characterization of $\mathscr{L}_{\omega}$-extension domains, see
 \cite{GM}. In \cite{L}, Lappalainen
generalized the characterization of \cite{GM}, and proved that $\Omega$ is an
$\mathscr{L}_{\omega}$-extension domain if and only if each pair of points
$x,y\in \Omega$ can be joined by a rectifiable curve $\gamma\subset \Omega$
satisfying

\be\label{eq1.0}
\int_{\gamma}\frac{\omega(d_{\Omega}(\zeta))}{d_{\Omega}(\zeta)}\,ds(\zeta)
\leq M\omega(|x-y|) \ee
with some fixed positive constant
$M=M(\Omega,\omega)$, where $ds$ stands for the arc length measure on
$\gamma$ and $d_{\Omega}(\zeta)$ is the Euclidean
 distance between $\zeta$ and the boundary of $\Omega$.  Furthermore, Lappalainen \cite[Theorem 4.12]{L} proved
that $\mathscr{L}_{\omega}$-extension domains  exist only for fast majorants
$\omega$. In particular, $\mathbb{B}^{n}$ is an $\mathscr{L}_{\omega}$-extension domain of $\mathbb{R}^n$ for the fast majorant $\omega$ (see \cite{Dy2}).



Hardy and Littlewood proved that if $\varphi\in\mathscr{L}_{\omega_{\beta}}(\mathbb{S}^{n-1})$, then the Poisson integral
of $\varphi$ belongs to $\mathscr{L}_{\omega_{\beta}}(\mathbb{B}^{n})$ (cf. \cite{AKM,Ch-18,GK,H-L-31,H-L}), where $\beta\in(0,1)$ and $\omega_{\beta}(t)=t^{\beta}$ for $t\geq0$.
In \cite[Lemma 4]{Dy1} and \cite[Lemma 1]{Dy2}, Dyakonov proved a more general Hardy and Littlewood type theorem as follows.

\begin{ThmA}\label{Lip-1} Let $n\geq2$ and $\omega$ be a slow majorant.
If $\varphi\in\mathscr{L}_{\omega}(\mathbb{S}^{n-1})$, then
 $P_{0}[\varphi]\in \mathscr{L}_{\omega}(\mathbb{B}^{n})$. 
\end{ThmA}

However, if $\alpha>0$, then the Hardy and Littlewood type theorem for  $P_{\alpha}[\varphi]$ is very different from   Theorem A for $\alpha=0$.
It is read as follows.




\begin{Thm}\label{thm-1} Let $n\geq2$, $\alpha>0$ and $\omega$ be a fast majorant.
If $\varphi\in\mathscr{L}_{\omega}(\mathbb{S}^{n-1})$, then
 $P_{\alpha}[\varphi]\in \mathscr{L}_{\omega}(\mathbb{B}^{n})$.
\end{Thm}


Combining Theorem A and Theorem \ref{thm-1} yields the following result.

\begin{Cor}
 Let $n\geq2$, $\omega$ be a majorant and $\varphi\in\mathscr{L}_{\omega}(\mathbb{S}^{n-1})$.
 \begin{enumerate}
\item[{\rm (I)}] If $\alpha=0$ and $\omega$ is a slow majorant, then $P_{0}[\varphi]\in \mathscr{L}_{\omega}(\mathbb{B}^{n})$;
\item[{\rm (II)}] If $\alpha>0$ and $\omega$ is a fast majorant, then
 $P_{\alpha}[\varphi]\in \mathscr{L}_{\omega}(\mathbb{B}^{n})$.
\end{enumerate}
\end{Cor}

For $a, b, c\in\mathbb{R}$ with $c\neq0, -1, -2, \ldots,$ the {\it
hypergeometric} function is defined by the power series
$$F(a,b;c;x)=\sum_{k=0}^{\infty}\frac{(a)_{k}(b)_{k}}{(c)_{k}}\frac{x^{k}}{k!},~|x|<1,$$
where $(a)_{0}=1$ and $(a)_{k}=a(a+1)\cdots(a+k-1)$ for $k=1, 2,
\ldots$ are the {\it Pochhammer} symbols. Obviously, for  $k=0, 1,
2, \ldots$, $(a)_{k}=\Gamma(a+k)/\Gamma(a)$ (cf. \cite{AAR}).
 In particular, for
 $a+b<c$, we have (see \cite[Theorem 2.2.2]{AAR}, \cite[p.234]{O-1})
 \be\label{eq-y}
F(a,b;c;1)=\lim_{x\rightarrow1^{-}}
F(a,b;c;x)=\frac{\Gamma(c)\Gamma(c-a-b)}{\Gamma(c-a)\Gamma(c-b)}<\infty.\ee
For $\alpha>-1$, we have (see \cite[(2.8)]{LP})
\beq\label{P-integral-1}
P_{\alpha}[1](x)&=&\int_{\mathbb{S}^{n-1}}\mathbf{P}_{\alpha}(x,\zeta)d\sigma(\zeta)\\ \nonumber
&=&C_{\alpha}F\left(-\frac{\alpha}{2},\frac{n}{2}-1-\frac{\alpha}{2};\frac{n}{2};|x|^{2}\right),
\eeq
which combined with \cite[Lemma 1.2]{O-1} implies that
\begin{equation}
\label{bound for Palpha1}
P_{\alpha}[1](x)\leq 1,
\quad x\in \mathbb{B}^n
\end{equation}
for $\alpha> n-2$.

Heinz in his classical paper \cite{HZ} proved the following result: if $f$ is a sense-preserving
harmonic diffeomorphism of $\mathbb{D}$ onto itself with $f(0)=0$, then
\be\label{Heinz}|f_{x}(z)|^{2}+|f_{y}(z)|^{2}\geq\frac{2}{\pi^{2}},~z\in\mathbb{D}.\ee
It uses the maximum principle of holomorphic functions and the following
sharp inequality
\be\label{Heinz-1}\liminf_{r\rightarrow1^{-}}\frac{|f(\zeta)-f(r\zeta)|}{1-r}\geq\frac{2}{\pi}\ee
to prove (\ref{Heinz}). Recently, Kalaj \cite[Theorem  2.5]{K} generalized (\ref{Heinz-1}) into the higher dimensional cases as follows.
If $f$ is a proper harmonic mapping of $\mathbb{B}^{n}$ onto itself such that $f(0)=0$, then
the following sharp inequality
\be\label{eq-Kalaj}\liminf_{r\rightarrow1^{-}}\frac{|f(\zeta)-f(r\zeta)|}{1-r}
\geq \frac{n!\big[1+n-(n-2)F\big(\frac{1}{2},1;\frac{n+3}{2};-1\big)\big]}{2^{\frac{3n}{2}}
\Gamma(\frac{n+1}{2})\Gamma(\frac{n+3}{2})},~|\zeta|=1\ee
holds. 
For $\alpha>0$ and $n=2$, from the proof of Chen and Vuorinen \cite[Theorem 2]{C-V},
we deduce that
if $\varphi$ is a continuous function in $\mathbb{S}^{1}$ with $\varphi(\mathbb{S}^{1})\subset \mathbb{S}^{1}$
such that $P_{\alpha}[\varphi](0)=0$, then $|P_{\alpha}[\varphi](x)|\leq P_{\alpha}[1](x)\leq 1$ by (\ref{bound for Palpha1})
and

\be\label{C-V-2015}\liminf_{r\rightarrow1^{-}}\frac{|P_{\alpha}[\varphi](\zeta)-P_{\alpha}[\varphi](r\zeta)|}{1-r}\geq\frac{\alpha}{2},~|\zeta|=1.\ee

We extend (\ref{Heinz-1}), (\ref{eq-Kalaj}) and (\ref{C-V-2015}) into  the following form.
Note that in the case $\alpha>n-2$,
we do not assume that $f(0)=0$.

\begin{Thm}\label{thm-2}
Let $\varphi=(\varphi_{1},\cdots,\varphi_{n})$ be a continuous mapping of $\mathbb{S}^{n-1}$ into $\mathbb{R}^{n}$ with $|\varphi|=1$
on $\mathbb{S}^{n-1}$, where   $n\geq2$.
For $\alpha\in[0,\infty)$, let $f=(P_{\alpha}[\varphi_{1}],\ldots,P_{\alpha}[\varphi_{n}])$.
 \begin{enumerate}
\item[{\rm (1)}] If $\alpha=0$ and $f(0)=0$, then the following sharp inequality
$$\liminf_{r\rightarrow1^{-}}\frac{|f(\zeta)-f(r\zeta)|}{1-r}\geq
\frac{n!\big[1+n-(n-2)F\big(\frac{1}{2},1;\frac{n+3}{2};-1\big)\big]}{2^{\frac{3n}{2}}
\Gamma(\frac{n+1}{2})\Gamma(\frac{n+3}{2})},~|\zeta|=1$$ holds. 
\item[{\rm (2)}] If $\alpha>n-2$, then
$$\liminf_{r\rightarrow1^{-}}\frac{|f(\zeta)-f(r\zeta)|}{1-r}\geq\frac{\alpha}{2}+1-\frac{n}{2},~|\zeta|=1.$$
\end{enumerate}
\end{Thm}

The proofs of Theorems  \ref{thm-1} and \ref{thm-2} will be presented in Sec. \ref{sec3}.



\section{The proof of the main results}\label{sec3}





\subsection{The proof of Theorem \ref{thm-1}}
Let $\zeta=(\zeta_{1},\ldots,\zeta_{n})\in\mathbb{S}^{n-1}$ and $x=|x|\xi=(x_{1},\ldots,x_{n})$, where $\xi\in\mathbb{S}^{n-1}$.
 For $j\in\{1,\cdots,n\}$, we have
\begin{equation}
\label{I1I2}
\left|\frac{\partial}{\partial x_{j}}P_{\alpha}[\varphi](x)\right|\leq I_1+I_2,
\end{equation}
where
\[
I_1=\int_{\mathbb{S}^{n-1}}\left|\frac{\partial}{\partial x_{j}}\mathbf{P}_{\alpha}(x,\zeta)\right||\varphi(\zeta)-\varphi(\xi)|d\sigma(\zeta)
\]
and
\[
I_2=|\varphi(\xi)|\left|\frac{\partial}{\partial x_{j}}P_{\alpha}[1](x)\right|.
\]
Since for $j\in\{1,\cdots,n\}$,
$$\frac{\partial}{\partial x_{j}}\mathbf{P}_{\alpha}(x,\zeta)=C_{\alpha}(1+\alpha)\frac{(1-|x|^{2})^{\alpha}}
{|x-\zeta|^{n+\alpha}}-C_{\alpha}(n+\alpha)\frac{(1-|x|^{2})^{1+\alpha}(x_{j}-\zeta_{j})}{|x-\zeta|^{n+\alpha+2}},$$
  we see that,
\beq\label{eq-r-1-0}
I_1&=&\int_{\mathbb{S}^{n-1}}\left|\frac{\partial}{\partial x_{j}}\mathbf{P}_{\alpha}(x,\zeta)\right||\varphi(\zeta)-\varphi(\xi)|d\sigma(\zeta)\\ \nonumber
&\leq&C_{\alpha}(1+\alpha)\int_{\mathbb{S}^{n-1}}\frac{|\varphi(\zeta)-\varphi(\xi)|(1-|x|^{2})^{\alpha}}{|x-\zeta|^{n+\alpha}}d\sigma(\zeta)\\ \nonumber
&&+C_{\alpha}(n+\alpha)\int_{\mathbb{S}^{n-1}}\frac{|\varphi(\zeta)-\varphi(\xi)|(1-|x|^{2})^{1+\alpha}}{|x-\zeta|^{n+\alpha+1}}d\sigma(\zeta).
\eeq
Note that
\beqq|\varphi(\zeta)-\varphi(\xi)|\leq\|\varphi\|_{\mathscr{L}_{\omega}(\mathbb{S}^{n-1}),s}\omega(|\zeta-\xi|)\eeqq
 and
\be\label{e-1-1}|\zeta-\xi|\leq|\xi-x|+|x-\zeta|=1-|x|+|x-\zeta|\leq2|x-\zeta|,\ee
which imply that
\be\label{eq-r-1-1}|\varphi(\zeta)-\varphi(\xi)|\leq
\|\varphi\|_{\mathscr{L}_{\omega}(\mathbb{S}^{n-1}),s}\omega(2|x-\zeta|)\leq2\|\varphi\|_{\mathscr{L}_{\omega}(\mathbb{S}^{n-1}),s}\omega(|x-\zeta|).\ee
It follows from (\ref{eq-r-1-0}) and (\ref{eq-r-1-1})   that

\beq\label{I1-2}
I_1
&\leq&2C_{\alpha}(1+\alpha)\|\varphi\|_{\mathscr{L}_{\omega}(\mathbb{S}^{n-1}),s}(1-|x|^{2})^{\alpha}
J_{1}\\ \nonumber
&&+2C_{\alpha}(n+\alpha)\|\varphi\|_{\mathscr{L}_{\omega}(\mathbb{S}^{n-1}),s}(1-|x|^{2})^{1+\alpha}
J_{2},
\eeq where $$J_{1}=\int_{\mathbb{S}^{n-1}}\frac{\omega(|x-\zeta|)}{|x-\zeta|^{n+\alpha}}d\sigma(\zeta)$$ and
$$J_{2}=\int_{\mathbb{S}^{n-1}}\frac{\omega(|x-\zeta|)}{|x-\zeta|^{n+\alpha+1}}d\sigma(\zeta).$$
We begin to estimate $J_{1}$ and $J_{2}$. Let
 $$E_{1}=\{\zeta\in\mathbb{S}^{n-1}:~|\zeta-\xi|\leq1-r\}~\mbox{and}~E_{2}=\{\zeta\in\mathbb{S}^{n-1}:~|\zeta-\xi|>1-r\},$$
where $r=|x|$.
By a rotation, we may assume that $\xi=(1,0,\cdots,0)$.
Then $\zeta=(\zeta_{1},\cdots,\zeta_{n})\in E_{1}$ has the following spherical coordinate transformation:
\beqq
\zeta_{1}&=&\cos\theta_{1},\\
\zeta_{2}&=&\sin\theta_{1}\cos\theta_{2},\\
&\vdots&\\
\zeta_{n-1}&=&\sin\theta_{1}\sin\theta_{2}\cdots\sin\theta_{n-2}\cos\theta_{n-1},\\
\zeta_{n}&=&\sin\theta_{1}\sin\theta_{2}\cdots\sin\theta_{n-2}\sin\theta_{n-1},
\eeqq
where $\theta_{1}\in[-2\arcsin\frac{1-r}{2},2\arcsin\frac{1-r}{2}],$ $\theta_{2},\cdots,\theta_{n-2}\in[0,\pi]$ and $\theta_{n-1}\in[0,2\pi]$.
Then
\beqq
\int_{E_{1}}d\sigma(\zeta)&=&2\int_{0}^{\kappa(r)}\sin^{n-2}\theta_{1}d\theta_{1}\int_{0}^{\pi}\sin^{n-3}\theta_{2}d\theta_{2}\cdots
\int_{0}^{\pi}\sin\theta_{n-2}d\theta_{n-2}\int_{0}^{2\pi}d\theta_{n-1}\\
&=&2\mathscr{H}_{n-2}\int_{0}^{\kappa(r)}\sin^{n-2}\theta_{1}d\theta_{1}\\
&\leq&2\mathscr{H}_{n-2}\kappa(r)\sin^{n-2}\kappa(r)\\
&\leq&2\mathscr{H}_{n-2}(1-r)^{n-2}\kappa(r),
\eeqq
where $\kappa(r)=2\arcsin\frac{1-r}{2}$ and $\mathscr{H}_{n-2}=2\pi^{(n-1)/2}/\Gamma((n-1)/2)$.
Since $$\lim_{r\rightarrow1^{-}}\frac{\kappa(r)}{1-r}=1,$$ there is a positive constant $M$ such that
\beqq\kappa(r)\leq M(1-r),\eeqq
which implies that
\be\label{eq-rt-1} \int_{E_{1}}d\sigma(\zeta)\leq 2M\mathscr{H}_{n-2}(1-r)^{n-1}.\ee
Since the function
\beqq
\frac{\omega(t)}{t^{n+\alpha}}=\frac{\omega(t)}{t}\cdot\frac{1}{t^{n+\alpha-1}}
\eeqq
is decreasing for $t\in(0,\infty)$,
by (\ref{eq-rt-1}), there is a positive constant $M$ such that
\be\label{eq-rt-2}
\int_{E_{1}}\frac{\omega(|x-\zeta|)}{|x-\zeta|^{n+\alpha}}d\sigma(\zeta)\leq\frac{\omega(1-r)}{(1-r)^{n+\alpha}}\int_{E_{1}}d\sigma(\zeta)
\leq 2M\mathscr{H}_{n-2}\frac{\omega(1-r)}{(1-r)^{1+\alpha}}.
\ee

Since, for $\zeta\in \mathbb{S}^{n-1}$,
\beqq
\frac{\omega(|x-\zeta|)}{|x-\zeta|}\leq\frac{\omega(1-r)}{1-r},
\eeqq by (\ref{e-1-1}) and \cite[p.317]{AKM}, we see that there is a positive   constant $M$ such that

\beq\label{eq-rt-3}
\int_{E_{2}}\frac{\omega(|x-\zeta|)}{|x-\zeta|^{n+\alpha}}d\sigma(\zeta)&\leq&\frac{\omega(1-r)}{1-r}\int_{E_{2}}\frac{1}{|x-\zeta|^{n+\alpha-1}}d\sigma(\zeta)\\ \nonumber
&\leq&2^{n+\alpha-1}\frac{\omega(1-r)}{1-r}\int_{E_{2}}\frac{1}{|\xi-\zeta|^{n+\alpha-1}}d\sigma(\zeta)\\ \nonumber
&\leq&M\frac{\omega(1-r)}{1-r}\int_{1-r}^{2}\rho^{-(n+\alpha-1)}\rho^{n-2}d\rho\\ \nonumber
&\leq&\frac{M}{\alpha}\frac{\omega(1-r)}{(1-r)^{1+\alpha}}.
\eeq
It follows  (\ref{eq-rt-2}) and (\ref{eq-rt-3}) that

\beq\label{eq-rt-4}
J_{1}&=&\int_{E_{1}}\frac{\omega(|x-\zeta|)}{|x-\zeta|^{n+\alpha}}d\sigma(\zeta)+\int_{E_{2}}\frac{\omega(|x-\zeta|)}{|x-\zeta|^{n+\alpha}}d\sigma(\zeta)\\ \nonumber
&\leq&
 \left(M\mathscr{H}_{n-2}+\frac{M}{\alpha}\right) \frac{\omega(1-r)}{(1-r)^{1+\alpha}},
 \eeq
where $M$ is a positive constant.

By using the similar reasoning as in the proof of $J_{1}$, we see that there a positive constant $M$ such that

\beq\label{eq-rt-5}
J_{2}&=&\int_{E_{1}}\frac{\omega(|x-\zeta|)}{|x-\zeta|^{n+\alpha+1}}d\sigma(\zeta)+\int_{E_{2}}\frac{\omega(|x-\zeta|)}{|x-\zeta|^{n+\alpha+1}}d\sigma(\zeta)\\ \nonumber
&\leq&
 M \frac{\omega(1-r)}{(1-r)^{2+\alpha}}.
 \eeq
By (\ref{I1-2}), (\ref{eq-rt-4}) and (\ref{eq-rt-5}), we conclude that there is a positive constant $M$ such that

\be\label{F-01} I_{1}\leq M \frac{\omega(1-r)}{1-r}.\ee

Next, we estimate $I_2$.
From (\ref{P-integral-1})
we have (see \cite[p.234]{O-1})
\[
\frac{\partial}{\partial x_{j}}P_{\alpha}[1](x)
=
\frac{2C_{\alpha}\left(-\frac{\alpha}{2}\right)\left(\frac{n}{2}-1-\frac{\alpha}{2}\right)}{\frac{n}{2}}x_j
F\left(-\frac{\alpha}{2}+1,\frac{n}{2}-\frac{\alpha}{2};\frac{n}{2}+1;|x|^{2}\right).
\]
Since $\alpha-\beta>-1$,
$P_{\alpha-\beta}[1](x)$ and $P_{1+\alpha-\beta}[1](x)$ are bounded on $\mathbb{B}^n$.
Also, since
\[
\frac{n}{2}+1-\left(-\frac{\alpha}{2}+1\right)-\left(\frac{n}{2}-\frac{\alpha}{2}\right)
=\alpha>0,
\]
$F\left(-\frac{\alpha}{2}+1,\frac{n}{2}-\frac{\alpha}{2};\frac{n}{2}+1;|x|^{2}\right)$ is bounded on $\mathbb{B}^n$.
Combining these results with (\ref{I1I2}), (\ref{F-01}) gives that
there is a positive constant $M>0$ such that
\be\label{eq-r-1-4}|\nabla P_{\alpha}[\varphi](x)|\leq M\frac{\omega(1-|x|)}{1-|x|}.
\ee
Since $\mathbb{B}^{n}$ is an $\mathscr{L}_{\omega}$-extension domain of $\mathbb{R}^n$ for the fast majorant $\omega$,
each pair of points
$x,y\in \mathbb{B}^{n}$ can be joined by a rectifiable curve $\gamma\subset\mathbb{B}^{n}$
satisfying
(\ref{eq1.0}).
Therefore, by (\ref{eq-r-1-4}), we conclude that
there is a positive constant $M>0$ such that
\beqq
|P_{\alpha}[\varphi](x)-P_{\alpha}[\varphi](y)|
&\leq&
\int_{\gamma}|\nabla P_{\alpha}[\varphi](u)|ds(u)
\\
&\leq &
M\int_{\gamma}\frac{\omega(1-|u|)}{1-|u|}ds(u)
\\
&\leq&
M\omega(|x-y|),
\eeqq
which implies that
$P_{\alpha}[\varphi]\in \mathscr{L}_{\omega}(\mathbb{B}^{n})$.
\qed

\subsection{The proof of Theorem \ref{thm-2}}
Since (1) easily follows from \cite[Theorem 2.5]{K}, we only need to prove (2).
By the assumptions, we have
\beqq
f(x)=\int_{\mathbb{S}^{n-1}}\mathbf{P}_{\alpha}(x,\zeta)\varphi(\zeta)d\sigma(\zeta),~x\in\mathbb{B}^{n}.
\eeqq Then, for $x\in\mathbb{B}^{n}$,
\beqq
|f(x)|&\leq& P_{\alpha}[1](x)=\int_{\mathbb{S}^{n-1}}\mathbf{P}_{\alpha}(x,\zeta)d\sigma(\zeta)\\
&=&C_{\alpha}F\left(-\frac{\alpha}{2},\frac{n}{2}-1-\frac{\alpha}{2};\frac{n}{2};|x|^{2}\right)\\
&=&C_{\alpha}\sum_{k=0}^{\infty}\frac{\left(-\frac{\alpha}{2}\right)_{k}\left(\frac{n}{2}-1-\frac{\alpha}{2}\right)_{k}}{(1)_{k}\left(\frac{n}{2}\right)_{k}}|x|^{2k},
\eeqq
which implies that
\be\label{eq-1}
\lim_{|x|\rightarrow1^{-}}P_{\alpha}[1](x)=1.
\ee
(Since $\alpha\geq 0$, (\ref{eq-1}) also follows from the fact that $P_{\alpha}[1]\in C(\overline{\mathbb{B}^n})$).
It follows from (\ref{eq-1}) and $\alpha>0$ that, for $r\in(0,1)$,

\beq\label{eq-2}\nonumber
\liminf_{r\rightarrow1^{-}}\frac{|f(\zeta)-f(r\zeta)|}{1-r}&\geq&\liminf_{r\rightarrow1^{-}}\frac{1-P_{\alpha}[1](r\zeta)}{1-r}=
\lim_{r\rightarrow1^{-}}\frac{1-P_{\alpha}[1](r\zeta)}{1-r}\\ \nonumber
&=&\lim_{r\rightarrow1^{-}}\frac{\partial P_{\alpha}[1](r\zeta)}{\partial r}\\ \nonumber
&=&2C_{\alpha}\lim_{r\rightarrow1^{-}}r\sum_{k=1}^{\infty}\frac{\left(-\frac{\alpha}{2}\right)_{k}
\left(\frac{n}{2}-1-\frac{\alpha}{2}\right)_{k}}{\left(\frac{n}{2}\right)_{k}}\frac{(r^{2})^{k-1}}{(k-1)!}\\ 
&=&\frac{4}{n}\left(-\frac{\alpha}{2}\right)\left(\frac{n}{2}-1-\frac{\alpha}{2}\right)
\mathcal{M}(n,\alpha),
\eeq where $$\mathcal{M}(n,\alpha)=C_{\alpha}F\left(-\frac{\alpha}{2}+1,\frac{n}{2}-\frac{\alpha}{2};\frac{n}{2}+1;1\right).$$
By (\ref{eq-y}), we have
\beqq
\mathcal{M}(n,\alpha)&=&\frac{\Gamma\left(\frac{n+\alpha}{2}\right)\Gamma\left(1+\frac{\alpha}{2}\right)}{\Gamma\left(\frac{n}{2}\right)\Gamma(1+\alpha)}
\frac{\Gamma\left(\frac{n}{2}+1\right)\Gamma\left(\alpha\right)}{\Gamma\left(\frac{n+\alpha}{2}\right)\Gamma\left(1+\frac{\alpha}{2}\right)}\\
&=&\frac{\frac{n}{2}\Gamma\left(\frac{n}{2}\right)\Gamma\left(\alpha\right)}{\alpha\Gamma\left(\frac{n}{2}\right)\Gamma\left(\alpha\right)}\\
&=&\frac{n}{2\alpha},
\eeqq which, together with (\ref{eq-2}), implies that
$$\liminf_{r\rightarrow1^{-}}\frac{|f(\zeta)-f(r\zeta)|}{1-r}\geq\frac{\alpha}{2}+1-\frac{n}{2},~|\zeta|=1.$$
The proof of this theorem is complete.
\qed



\section*{Statements and Declarations}

\subsection*{Competing interests}
There are no competing interests.

\subsection*{Data availability}
Data sharing not applicable to this article as no datasets were generated or analysed during the current study.

\bigskip
\section*{Acknowledgments}
The research of the first author was partly supported by the National Science
Foundation of China (No. 12071116), the Hunan Provincial
Natural Science Foundation of China (No. 2022JJ10001), the Key
Projects of Hunan Provincial Department of Education (grant no.
21A0429),
 the Double First-Class University Project of Hunan Province
(Xiangjiaotong [2018]469),  the Science and Technology Plan Project
of Hunan Province (No. 2016TP1020),  and the Discipline Special Research
Projects of Hengyang Normal University (No. XKZX21002); The research of
the second author was partly supported by JSPS KAKENHI Grant Number
JP22K03363; The third author was supported by the extracurricular academic and
scientific works competition project of Hengyang Normal University (No.33).


\end{document}